\documentclass[12pt]{article}

\usepackage{amsmath,amssymb,amsthm,amscd}
\usepackage[colorlinks, linkcolor=blue, anchorcolor=gree, citecolor=red]{hyperref}
\usepackage[numbers,sort&compress]{natbib}
\usepackage{graphicx}

\begin{document}

\newcommand{\Cyc}{{\rm{Cyc}}}

\newtheorem{thm}{Theorem}[section]
\newtheorem{pro}[thm]{Proposition}
\newtheorem{lem}[thm]{Lemma}
\newtheorem{cor}[thm]{Corollary}
\theoremstyle{definition}
\newtheorem{ex}[thm]{Example}
\newtheorem{remark}[thm]{Remark}
\newcommand{\bth}{\begin{thm}}
\renewcommand{\eth}{\end{thm}}
\newcommand{\bex}{\begin{examp}}
\newcommand{\eex}{\end{examp}}
\newcommand{\bre}{\begin{remark}}
\newcommand{\ere}{\end{remark}}

\newcommand{\bal}{\begin{aligned}}
\newcommand{\eal}{\end{aligned}}
\newcommand{\beq}{\begin{equation}}
\newcommand{\eeq}{\end{equation}}
\newcommand{\ben}{\begin{equation*}}
\newcommand{\een}{\end{equation*}}

\newcommand{\bpf}{\begin{proof}}
\newcommand{\epf}{\end{proof}}
\renewcommand{\thefootnote}{}

\def\beql#1{\begin{equation}\label{#1}}
\title{\Large\bf Finite groups with star-free noncyclic graphs}

\author{{\sc Xuanlong Ma$^{1,3}$, Gary L. Walls$^2$, Kaishun Wang$^3$
}\\[15pt]
{\small\em $^1$College of Mathematics and Information Science, Guangxi University,
Nanning, 530004, China}\\
{\small\em $^2$Department of Mathematics, Southeastern Louisiana University, Hammond, LA 70402, USA}\\
{\small\em $^3$Sch. Math. Sci. {\rm \&} Lab. Math. Com. Sys., Beijing Normal University, Beijing, 100875, China}\\
}

 \date{}

\maketitle

\begin{abstract}
For a finite noncyclic group $G$, let $\Cyc(G)$ be a set of elements $a$ of $G$ such that $\langle a,b\rangle$ is cyclic for each  $b$ of $G$.  The noncyclic graph of $G$ is a graph with the vertex set $G\setminus \Cyc(G)$, having an edge between two
distinct vertices $x$ and $y$ if $\langle x, y\rangle$ is not cyclic. In this paper,
we classify all finite noncyclic groups whose noncyclic graphs are $K_{1,n}$-free, where $3\leq n\leq 6$.
\end{abstract}



\footnote{E-mail addresses: xuanlma@mail.bnu.edu.cn (X. Ma), gary.walls@selu.edu (G.L. Walls) wangks@bnu.edu.cn (K. Wang).}
\section{Introduction}

All groups considered in this paper are finite.
Let $G$ be a noncyclic group. The {\em cyclicizer} $\Cyc(G)$ of $G$ is the set
$$\{a\in G: \langle a,b\rangle\ \text{is cyclic for each }  b\in G\},$$
which is a normal cyclic subgroup of $G$ (see \cite{OP92}).
The noncyclic graph $\Gamma_G$ of $G$
is the graph whose vertex set is $G\setminus \Cyc(G)$, and two distinct vertices being adjacent if they do not generate a cyclic subgroup.
In 2007, Abdollahi and  Hassanabadi \cite{AH07} introduced the concept of noncyclic graphs and
established basic graph theoretical properties. Aalipour et. al  \cite{acam} discovered the relationship between the complement graph of a noncyclic graph and two well-studied graphs--power graphs and commuting graphs. For the latter two graphs, see \cite{KQ,CGS,BF}.
In \cite{AH09}, Abdollahi and  Hassanabadi investigated the clique number  of a noncyclic graph. The full automorphism group of a noncyclic graph was characterized in
\cite{MLW1} and the noncyclic graphs of genus at most four was classified in \cite{MLW2}.

A graph is said to be {\em $\Gamma$-free} if
it has no induced subgraphs isomorphic to $\Gamma$.
Forbidden graph characterization
appears in many contexts; for instance, forbidden subgraph problem (Tur\'{a}n-type problem), or extremal graph theory where
lower and upper bounds can be obtained for various numerical invariants of the corresponding
graphs.
Some graphs obtained from groups with small forbidden induced subgraphs have been studied in the literature. For example, Doostabadi et al. \cite{DoE} studied the
power graphs of finite groups with  $K_{1,3}$ or $C_4$-free.
Akhlaghi and Tong-Viet \cite{AkT} studied the finite groups with $K_4$-free prime graphs. Rajkumar and Devi \cite{RD15} classified the finite groups with $K_4$ or $K_5$-free intersection graphs of subgroups.
In \cite{DaN}, Das and Nongsiang classified
$K_3$-free commuting graphs of finite non-abelian groups.

In this paper, we study noncyclic graphs of  finite groups.
In Section \ref{sec2} we classify all finite groups $G$ with a unique involution and   $\pi_e(G)=\{2,3,4,6\}$.
In Section \ref{sec3},
we classify all finite noncyclic groups whose noncyclic graphs are $K_{1,n}$-free, where $3\leq n\leq 6.$

\section{A result about finite groups}\label{sec2}

An element of order $2$ in a
group is called an {\em involution}.
The set of natural numbers consisting of orders of non-identity elements of $G$ is denoted by $\pi_e(G)$. The {\em exponent} of $G$ is
the least common multiple of the orders of the elements of $G$.
We denote the cyclic group of order $n$ and the quaternion group of
order $8$ by $\mathbb{Z}_n$ and $Q_8$, respectively.
Also $\mathbb{Z}_n^m$ is used for the $m$-fold direct product of the cyclic group $\mathbb{Z}_n$ with itself.

In this section we prove the following result about finite groups, which will be used to classify finite groups with $K_{1,5}$-free noncyclic graphs.

\begin{thm}\label{afgt}
Let $G$ be a finite group having a unique involution and that $\pi_e(G)=\{2,3,4,6\}$.
Then, either $G\cong SL(2,3)$ or $G\cong \mathbb{Z}_3^n\rtimes \mathbb{Z}_4$, where
$\mathbb{Z}_4$ acts on $\mathbb{Z}_3^n$ by inversion.
\end{thm}

Let $G$ be a finite group and $p$ be a prime number dividing $|G|$. Denote
by ${\rm Syl}_p(G)$ and $O_p(G)$ the set of all Sylow $p$-subgroups, and
the largest normal $p$-subgroup of $G$, respectively. Note that
$O_p(G)=\bigcap_{P\in {\rm Syl}_p(G)} P$.
Let $n_p=|{\rm Syl}_p(G)|$
and $P\in {\rm Syl}_p(G)$.
Recall that $n_p=|G: N_G(P)|\equiv 1 \pmod{p}$
and $n_p$ is a divisor of $|G:P|$.

\begin{lem}\label{opg} Let $G$ be a finite group and suppose that $n_p=p+1$ for some prime number $p$. Then for any two distinct $P_i,P_j\in {\rm Syl}_p(G)$, $P_i\cap P_j=O_p(G)$.
\end{lem}
\bpf
Let $m=n_p$ and $L={\rm Syl}_p(G)=\{P_1,P_2,\cdots,P_m\}$.
Now in order to prove the required result, it suffices to prove $P_1\cap P_2=O_p(G)$.

Let $R=P_1\cap P_2$ and let
$R$ act on $L$ by conjugation. Since
$(R\cap N_G(P_i))P_i=P_i(R\cap N_G(P_i))$ for all $i$, one has $R\cap N_G(P_i)\subseteq P_i$.
It follows that $R_{P_i}=R\cap N_G(P_i)= R\cap P_i$, where $R_{P_i}$ is the stabilizer of $P_i$
in $R$. Note that $R=P_1\cap P_2$. Then ${\rm Orbit}_R(P_1)|=|{\rm Orbit}_R(P_2)|=1$, where ${\rm Orbit}_R(P_i)$ is the $R$-orbit that contains $P_i$. Note that every $R$-orbit has length $1$ or $p$. Since $|L|=p+1$, any $R$-orbit has length $1$. This implies that
$R=R_{P_i}$  for all $i$. It follows that for each $i\ge 3$, $P_1\cap P_2=P_1\cap P_2 \cap P_i$,
namely, $P_1\cap P_2\subseteq P_i$. Thus, $P_1\cap P_2 \subseteq \bigcap_{i=3}^m P_i$ and
so $P_1\cap P_2=\bigcap_{P\in {\rm Syl}_p(G)} P=O_p(G)$, as desired.
\epf

We now note that for any prime number $p$, a $p$-group with a unique subgroup of order $p$ is either a cyclic group or a generalized quaternion group (cf. \citep[Theorem 5.4.10 (ii)]{Gor}).

\bigskip

{\noindent \em Proof of Theorem \ref{afgt}.}
Suppose that $|G|=2^t\cdot 3^n$ for some $t\ge 1,n\ge 1$.
Let $Q$ and $P$ be a Sylow $2$-subgroup and a Sylow $3$-subgroup of $G$, respectively.
Since $G$ has a unique involution and $8\notin \pi_e(G)$, we know that $Q\in \{Q_8, \mathbb{Z}_4\}$. Since $G$ has no elements of order $9$,
one has that $P$ has exponent $3$. Denote by $x$  the unique involution of $G$. Then $x\in Z(G)$, the center of $G$.

\medskip
\noindent {\bf Case 1.} $Q=Q_8$.
\medskip

Let $\langle a\rangle$, $\langle b\rangle$, and $\langle c\rangle$ be the three cyclic subgroups of $Q$ of order $4$, and that $ab=c$. Then $a^2=b^2=c^2=x$.
Since in this case $|G|=8\cdot 3^n$, we have that $n_3$ is a divisor of $8$. This implies that $n_3=1$ or $4$.

Suppose that $O_3(G)\ne 1$. Then let $a,b,c$ act on $O_3(G)$ by conjugation.
Neither of them can fix any non-identity elements of $O_3(G)$, since $G$ has no elements of order $12$.  Thus, $a,b,c$ act as fixed-point-free automorphisms of $O_3(G)$. Since $a^2=b^2=c^2\in Z(G)$,
one has that $a,b,c$ act as fixed-point-free automorphisms of order $2$. Now by Burnside's result (see \cite{Bu55}), we know that $O_3(G)$ is abelian and for each non-trivial element $g\in O_3(G)$,
we have that
$g^a=g^{-1}=g^b=g^c$. It then follows that $g^{c}=g^{ab}=g$, and hence $|abg|=12$, a contradiction. Therefore, we get $O_3(G)= 1$.

Thus, now we know $n_3=4$. By Lemma \ref{opg} we have that for any two
distinct $P_i,P_j\in {\rm Syl}_3(G)$, $P_i\cap P_j=1$. It follows that
the number of elements of order $3$ is $4(3^n-1)$. Also, since every element of order $3$
and $x$ can generate a cyclic subgroup of order $6$, the number of elements of order $6$ is also
$4(3^n-1)$. Now all that remains is to count the number of elements of order $4$.

Let $w$ be an element of order $4$ in $G$. Then, there is a $Q_1\in {\rm Syl}_2(G)$
so that $w\in Q_1$. Note that $Q_1\cong Q_8$. It follows that $Q_1\subseteq N_G(\langle w\rangle)$. If there exists an element $y$ of order $3$ such that $\langle w\rangle^y=\langle w\rangle$, then $\langle w\rangle$ is normal in $\langle w\rangle\langle y\rangle$ and $|\langle w\rangle\langle y\rangle|=12$, which implies that $\langle w\rangle\langle y\rangle\cong \mathbb{Z}_{12}$, a contradiction.
It follows that $Q_1= N_G(\langle w\rangle)$.
Thus, every element of order $4$ is contained in a unique Sylow $2$-subgroup of $G$.
It means that the number of elements of order $4$ is $6n_2$.

Suppose that $N_G(Q)=Q$. Then $n_2=3^n$. Counting all the elements of $G$ gives that
$$
8\cdot 3^n=6\cdot 3^n+8(3^n-1)+2.
$$
This implies that $3^n=1$, contrary to the order of $G$.
Thus, we obtain that $Q\subset N_G(Q)$.

Suppose that $|P\cap N_G(Q)|\ge 9$. Then there exist $w_1,w_2\in P\cap N_G(Q)$ so that
$\langle w_1,w_2\rangle$ is an abelian group of order $9$.
Now both $w_1$ and $w_2$ act on $Q$ by conjugation.
If follows that both $w_1$ and $w_2$ act as $3$-cycles on $\{\langle a\rangle, \langle b\rangle, \langle c\rangle\}$. But then there is an element $u$ of order $3$ that fixes some
cyclic subgroup $\langle v\rangle$ of order $4$,
where $u=w_1w_2^i$ for some integer $i$. It follows that there exist an
element of order $12$ in $\langle u\rangle\langle v\rangle$, a contradiction. Thus, we get
$|P\cap N_G(Q)|\le 3$.

Note that $N_G(Q)=G\cap N_G(Q)=Q(P\cap N_G(Q))$. Since $Q\subset N_G(Q)$, one has that
$P\cap N_G(Q)$ is a subgroup of order $3$ and $|N_G(Q)|=24$. This forces that $n_2=3^{n-1}$.
Now as above we get
$$
8\cdot 3^n=6\cdot 3^{n-1}+8(3^n-1)+2,
$$
which implies $n=1$ and so $|G|=24$.
Note that in this case $Q$ is normal in $G$. It is easy to see that $G\cong SL(2,3)$.

\medskip
\noindent {\bf Case 2.} $Q=\mathbb{Z}_4$.
\medskip

Let $Q=\langle y\rangle$. Since $\langle x\rangle P\subseteq N_G(P)$,
one has that $|G:N_G(P)|\ne 4$. Note that $n_3$ is a divisor of $4$. Then
$n_3=1$ and so $P$ is normal in $G$. Now as above $y$ acts as a fixed-point-free automorphism
of order $2$ on $P$ by conjugation. By Burnside's result, $P$ is abelian and so $P\cong \mathbb{Z}_3^n$ for some $n$, and for all $w\in P$ we have $w^y=w^{-1}$. It follows that
$G\cong \mathbb{Z}_3^n\rtimes \mathbb{Z}_4$, as desired.
$\qed$
\section{Main results}\label{sec3}
In this section we classify all finite groups with $K_{1,n}$-free noncyclic graphs, where 
$3\le n \le 6$.

We always use $G$ to denote a finite noncyclic group with the identity element $e$.  Euler's totient function is denoted by $\phi$.
A proper cyclic subgroup $\langle x\rangle$
is said to be {\em maximal} in $G$ if $\langle x\rangle \subseteq \langle y\rangle$
implies that $\langle x\rangle=\langle y\rangle$ or $G=\langle y\rangle$, where $y$ is an element of $G$. 
We first begin with the following two lemmas which will be used frequently in the sequel.

\begin{lem}\label{k1n} Suppose that $\langle g\rangle$ is
a maximal cyclic subgroup of $G$. Then
$\Gamma_G$ has an induced subgraph isomorphic to $K_{1,\phi(|g|)}$.
\end{lem}
\bpf
Let $n=\phi(|g|)$ and let $\{g_1,g_2,\cdots,g_n\}$ be all generators of $\langle g\rangle$.
Note that $G$ is noncyclic. Pick an element $a$ in $G\setminus\langle g\rangle$.
Since $\langle g\rangle$ is maximal cyclic, one has that
$\langle a,g_i\rangle$ is not cyclic for each $i\in \{1,2,\ldots,n\}$. This implies that
$\{g_1,g_2,\cdots,g_n,a\}$ induces a subgraph isomorphic to $K_{1,n}$.
\epf

For a
graph $\Gamma$, we denote the sets of the vertices and the edges of ¦£ by $V(\Gamma)$ and $E(\Gamma)$, respectively.
An {\em independent set} is
a set of vertices in a graph, no two of which are adjacent; that is, a set whose induced
subgraph is null. The {\em independence number} of a graph $\Gamma$ is the cardinality of the largest independent set and is denoted by $\alpha(\Gamma)$.
The following result follows from \cite[Proposition 4.6]{AH07}.

\begin{lem}\label{indp} $\alpha(\Gamma_G)=\max\{|g|: g\in G\}-|\Cyc(G)|$.
\end{lem}

$\Gamma_G$ is complete if and only if $G$
is an elementary abelian $2$-group (see \cite[Proposition 3.1]{AH07}). So, we first note that
$\Gamma_G$ is $K_{1,2}$-free if and only if $G$ is an elementary abelian $2$-group.

Now we classify the finite groups whose noncyclic graphs are claw-free ($K_{1,3}$-free).

\begin{thm}\label{thk13}
$\Gamma_{G}$ is claw-free if and only if
$G$ is isomorphic to one of the following groups:

$(a)$ $Q_8$;

$(b)$ $\mathbb{Z}_2^n$, $n\ge 2$;

$(c)$  A noncyclic $3$-group of exponent $3$;

$(d)$ A noncyclic group $G$ with $\pi_e(G)=\{2,3\}$.
\end{thm}
\bpf
Since $\Gamma_{Q_8}\cong K_{2,2,2,2}$, $\Gamma_{Q_8}$ is claw-free. By Lemma \ref{indp} we
see that the independence number of the noncyclic graph of every group
in $(b)$--$(d)$ is at most $2$, and so each of the noncyclic graphs is claw-free.

Now we suppose that $\Gamma_{G}$ is claw-free.
Then, it follows from Lemma \ref{k1n} that $\phi(|g|)\le 2$ for
any maximal cyclic subgroup $\langle g\rangle$ of $G$.
This implies that every cyclic subgroup of $G$ has at most two generators.
Thus, $\pi_e(G)\subseteq\{2,3,4,6\}$.

Suppose that $G$ has an element $a$ of order $6$. Then according to Lemma \ref{k1n},
$\langle a\rangle$ is maximal. Note that $G$ is noncyclic.
Pick an element $x$ in $G\setminus \langle a\rangle$. If $|x|=2$,
then $\{x,a,a^3,a^5\}$ induces a subgraph isomorphic to $K_{1,3}$, a contradiction.
If the order of $|x|$ is $3$ or $4$, then $\{x,a,a^2,a^5\}$ also induces a subgraph isomorphic to $K_{1,3}$.
If $|x|=6$, then $\langle x\rangle \cap \langle a\rangle=2$ or $3$, and so in $G\setminus \langle a\rangle$, $G$ has an element of order $2$ or $3$,  a contradiction.

Now we know that $\pi_e(G)\subseteq\{2,3,4\}$. Suppose that there exists an element
$g$ of order $4$ in $G$.
If $x\in G\setminus \langle g\rangle$ and $|x|=2$ or $3$, then $\{x,g,g^2,g^3\}$ induces a subgraph isomorphic to $K_{1,3}$, a contradiction. Consequently, in this case we have
that $G$ has a unique involution and $\pi_e(G)=\{2,4\}$.
It follows that $G$ is isomorphic to $Q_8$.

Thus, we may assume that $\pi_e(G)\subseteq\{2,3\}$. If $\pi_e(G)=\{2\}$, then
$G$ is an elementary abelian $2$-group. If $\pi_e(G)=\{3\}$, then
$G$ is a $3$-group of exponent $3$.
\epf

\begin{thm}\label{thk14}
$\Gamma_{G}$ is $K_{1,4}$-free if and only if
$G$ is isomorphic to one of the following groups:

$(a)$ A noncyclic group $G$ with $\pi_e(G)\subseteq\{2,3,4\}$;

$(b)$ $\mathbb{Z}_6\times \mathbb{Z}_2^m,~~m\ge 1$.
\end{thm}
\bpf
If $G$ is isomorphic to a noncyclic group with $\pi_e(G)\subseteq\{2,3,4\}$,
then $\alpha(\Gamma_G)\le 3$ by Lemma \ref{indp}, and so $\Gamma_{G}$ is $K_{1,4}$-free.
If $G\cong \mathbb{Z}_6\times \mathbb{Z}_2^m$ for some $m\ge 1$, then
$|\Cyc(G)|=3$ and so $\Gamma_G$ is a complete multipartite graph whose each partite set has size $3$, which implies that in this case $\Gamma_{G}$ is also $K_{1,4}$-free.

Now suppose that $\Gamma_{G}$ is $K_{1,4}$-free. Note that $\phi(n)$ is even for
any integer $n\ge 3$. By Lemma \ref{k1n} we see that each cyclic subgroup
of $G$ has at most two generators.
It follows that $\pi_e(G)\subseteq\{2,3,4,6\}$. In order to get desired result,
we now suppose that $G$ has an element $g$ of order $6$.
Clearly, $\langle g\rangle$ is maximal cyclic.
If there exists an element $a$ in $G\setminus \langle g\rangle$ such that $|a|=3$ or $4$,
then $\{a,g,g^2,g^4,g^5\}$ induces a subgraph isomorphic to $K_{1,4}$, a contradiction.
This means that $G$ has a unique subgroup of order $3$ and $\pi_e(G)=\{2,3,6\}$.

Let $P$ and $Q$ be a Sylow $2$-subgroup and a Sylow $3$-subgroup, respectively.
Then $G=P\ltimes Q$, $P$ is an elementary abelian $2$-group
of order great than $2$ and $Q\cong \mathbb{Z}_3$. Pick an involution $u$ in $G$.
If $\langle u,g^2\rangle$ is noncyclic, then $\{u,g,g^2,g^4,g^5\}$ induces a subgraph isomorphic to $K_{1,4}$. Thus, we have that every element of $P$ and every element of $Q$ commute.
It follows that
$$G= P\times Q\cong \mathbb{Z}_6\times \mathbb{Z}_2^m,~~m\ge 1,$$
as required.
\epf

\begin{thm}\label{thk15}
$\Gamma_{G}$ is $K_{1,5}$-free if and only if
$G$ is isomorphic to one of the following groups:

$(a)$ A noncyclic group $G$ with $\pi_e(G)\subseteq\{2,3,4,5\}$;

$(b)$ $\mathbb{Z}_6\times \mathbb{Z}_2^m,~~m\ge 1$;

$(c)$ $\mathbb{Z}_2\times Q$, where $Q$ is a noncyclic $3$-group of exponent $3$;

$(d)$ The special linear group $SL(2,3)$;

$(e)$ $\mathbb{Z}_3^n\rtimes \mathbb{Z}_4$, where
$\mathbb{Z}_4$ acts on $\mathbb{Z}_3^n$ by inversion and $n\ge 1$.
\end{thm}
\bpf
If $G$ is isomorphic to one group in $(a)$,
then $\alpha(\Gamma_G)\le 4$ by Lemma \ref{indp}, and so $\Gamma_{G}$ is $K_{1,5}$-free.
By Theorem \ref{thk14}, $\Gamma_{\mathbb{Z}_6\times \mathbb{Z}_2^m}$ is $K_{1,5}$-free,
where $m\ge 1$. If $G\cong \mathbb{Z}_2\times Q$ for some noncyclic $3$-group $Q$ of exponent $3$, then $|\Cyc(G)|=2$ and so $\Gamma_G$ is a complete multipartite graph whose each partite set has size $4$, which implies that in this case $\Gamma_{G}$ is also $K_{1,5}$-free.
Now it is easy to see that $\Gamma_G$ is a complete multipartite graph whose maximal partite set has size $4$ if $G$ is isomorphic to
$SL(2,3)$ or one group in $(f)$. This means that
$\Gamma_{G}$ is $K_{1,5}$-free if $G$ is one group of $(d)$ and $(e)$.

Now suppose that $\Gamma_{G}$ is $K_{1,5}$-free. It follows from Lemma \ref{k1n} that
$\pi_e(G)\subseteq \{2,3,4,5,6,8,10,12\}$.

Suppose that $g\in G$ with $|g|=12$. If there exists an element $a$ with order
$5$ or $8$, then $\{a,g,g^2,g^5,g^7,g^{11}\}$ induces a subgraph isomorphic
to $K_{1,5}$, a contradiction. This implies that $\pi_e(G)\subseteq \{2,3,4,6,12\}$.
If $G$ has an element $b$ of $G\setminus \langle g\rangle$ with $|b|<12$,
then $\{b,g,g^5,g^7,g^{11},g^t\}$ induces a subgraph isomorphic
to $K_{1,5}$, where $|g^t|=|b|$. Thus, in this case one has $G\cong \mathbb{Z}_{12}$,
a contradiction. This means that $G$ has no elements of order $12$.
Similarly, we can get $10\notin \pi_e(G)$. If $G$ has an element of order $8$, then
a similar argument implies that $G$ is a $2$-group and it has a unique involution and a unique
cyclic subgroup of order $4$,
which implies that $G$ is a generalized quaternion group which has precisely two elements of
order $4$, a contradiction. Thus, now we have $\pi_e(G)\subseteq \{2,3,4,5,6\}$.

Now in order to get desired result,
we suppose that $G$ has an element $h$ of order $6$. Then it is easy to see that
$5\notin \pi_e(G)$.

\medskip
\noindent {\bf Case 1.}  $G$ has two distinct cyclic subgroups of order $6$.
\medskip

We assume that $H_1=\langle h\rangle$ and $H_2=\langle h_2\rangle$ are two distinct cyclic subgroups of order $6$. In order to avoid that $\{h_2,h,h^2,h^3,h^4,h^5\}$ induces $K_{1,5}$,
we may assume that $|H_1\cap H_2|\ge 2$. In fact, any two distinct cyclic subgroups of
order $6$ have non-trivial intersection.

\medskip
\noindent {\bf Subcase 1.1.}  There exist two distinct cyclic subgroups of
order $6$ such that their intersection has order $3$.
\medskip

Without loss of generality, let $|H_1\cap H_2|= 3$.
If $G$ has an element $x$ of order $4$, then without loss of generality, we may
assume that $\langle h^3,x\rangle$ is not cyclic, which implies that
$\{x,h,h^2,h^3,h^4,h^5\}$ induces a subgraph isomorphic to $K_{1,5}$, a contradiction.
Hence, in this subcase $\pi_e(G)\subseteq \{2,3,6\}$. Since every generator of any maximal cyclic subgroup of order $2$ or $3$ is adjacent to each of $\langle h\rangle\setminus \{e\}$,
every cyclic subgroup of order $2$ or $3$ is not maximal.
If $\langle y\rangle\ne \langle h^2\rangle$ is a subgroup of order $3$, and let
$\langle y\rangle\subseteq \langle h_3\rangle$ with $|h_3|=6$, then
$|\langle h_3\rangle\cap H_i|=2$ for $i=1,2$ and so
$h_3^3=h^3=h_2^3$, which is impossible as $H_1\ne H_2$.
This implies that $G$ has a unique subgroup of order $3$. Now we know
that $G\cong \mathbb{Z}_2^m\ltimes \mathbb{Z}_3$ for some integer $m\ge 2$.
Pick an involution $u$ in $G$. Then since $\langle u\rangle$ is not maximal,
there exists an element $h'$ of order $6$ such
that $\langle u\rangle\subseteq \langle h'\rangle$. Note that the uniqueness of the
subgroup of order $3$. Then $\langle h'\rangle\cap \langle h\rangle=\langle h^2\rangle$.
It follows that $\langle u,h^2\rangle$ is cyclic. Namely, every involution of $G$ and
$h^2$ commute. This implies that $G\cong \mathbb{Z}_2^m\times \mathbb{Z}_3$ for some integer $m\ge 2$, as desired.

\medskip
\noindent {\bf Subcase 1.2.}  The intersection of each two distinct cyclic subgroups of
order $6$ has order $2$.
\medskip

In this case we first claim that $G$ has a unique involution. Assume, to the contrary, that
$u$ is an involution of $G$ such that $u\ne h^3$.
Then $\langle u, h^2\rangle$ is not cyclic, since there are no two cyclic subgroups of
order $6$ such that their intersection has order $3$. This implies that
$\{u,h,h^2,h^3,h^4,h^5\}$ induces a subgraph isomorphic to $K_{1,5}$, a contradiction.
Thus, our claim is valid.

Now note that $\pi_e(G)\subseteq \{2,3,4,6\}$. If $4\notin \pi_e(G)$, then
$G\cong \mathbb{Z}_2\times Q$, where $Q$ is a noncyclic $3$-group of exponent $3$.
Thus, we may assume that $\pi_e(G)= \{2,3,4,6\}$.
Note that $\mathbb{Z}_3\rtimes \mathbb{Z}_4$ has a unique cyclic subgroup of order $6$,
where $\mathbb{Z}_4$ acts on $\mathbb{Z}_3$ by inversion.
By Theorem \ref{afgt} we see that that $G\cong SL(2,3)$ or $G\cong \mathbb{Z}_3^n\rtimes \mathbb{Z}_4$,
where $\mathbb{Z}_4$ acts on $\mathbb{Z}_3^n$ by inversion, and $n\ge 2$,
as required.

\medskip
\noindent {\bf Case 2.}  $G$ has a unique cyclic subgroups of order $6$.
\medskip

First we see that $\langle h\rangle$ is a normal subgroup of $G$.
Note that $\pi_e(G)\subseteq\{2,3,4,6\}$.
If $G$ has an element $x$ in $G\setminus \langle h\rangle$ such that
$x\in C_G(h)$, the centralizer of $h$ in $G$,
then $G$ has a subgroup isomorphic to $\mathbb{Z}_2\times \mathbb{Z}_6$ or $\mathbb{Z}_3\times \mathbb{Z}_6$, which contradicts the fact that $G$ has precisely two elements of order $6$. This implies that $C_G(h)=\langle h\rangle$. So $G/\langle h\rangle$ is isomorphic to a subgroup of $\mathbb{Z}_2$, and hence $|G|=6$ or $12$.
By verifying we get that $G\cong \mathbb{Z}_3\rtimes \mathbb{Z}_4$
where $\mathbb{Z}_4$ acts on $\mathbb{Z}_3$ by inversion,
as desired.
\epf

\begin{thm}\label{thk16}
$\Gamma_{G}$ is $K_{1,6}$-free if and only if
$G$ is isomorphic to one of the following groups:

$(a)$ A noncyclic group $G$ with $\pi_e(G)\subseteq\{2,3,4,5,6\}$;

$(b)$ $\mathbb{Z}_{10}\times \mathbb{Z}_2^m,~~m\ge 1$.
\end{thm}
\bpf
If $G$ is isomorphic to one group in $(a)$, then it follows from Lemma \ref{indp} that
$\alpha(\Gamma_G)\le 5$, and so $\Gamma_{G}$ is $K_{1,6}$-free.
If $G\cong \mathbb{Z}_{10}\times \mathbb{Z}_2^m$ for some $m\ge 1$, then $|\Cyc(G)|=5$ and so $\Gamma_G$ is a complete multipartite graph whose each partite set has size $5$, which implies $\Gamma_{G}$ is $K_{1,6}$-free.

Now suppose that $\Gamma_{G}$ is $K_{1,6}$-free. Since $\phi(n)$ is even for $n\ge 3$,
by Lemma \ref{k1n} one has $\phi(|g|)\le 4$ for any $g\in G$.
It follows that $\pi_e(G)\subseteq \{2,3,4,5,6,8,10,12\}$.
An argument similar to the one
used to third paragraph of the proof of Theorem \ref{thk15} shows that $8,12\notin \pi_e(G)$.
Consequently, we have $\pi_e(G)\subseteq \{2,3,4,5,6,10\}$.

Now in order to get desired result,
we suppose that $G$ has an element $h$ of order $10$. Then it is easy to see that
$\pi_e(G)= \{2,5,10\}$ and $G$ has a unique subgroup of order $5$.
Thus, we may assume that $G\cong P\ltimes \mathbb{Z}_5$, where $P$ is an elementary
abelian $2$-group of order at least $4$. Pick any involution $u$ in $P$.
If $\langle u, h^2\rangle$ is not cyclic, then $\{u,h^2,h^4,h^6,h^8,h^5,h\}$ induces a subgraph
isomorphic to $K_{1,6}$, a contradiction. Thus, every element in $P$ and $h^2$ commute.
It follows that $G\cong P\times \mathbb{Z}_5$, that is, $G\cong \mathbb{Z}_{10}\times \mathbb{Z}_2^m$ for some $m\ge 1$, as desired.
\epf

\section*{Acknowledgement}
K. Wang's research is supported by National Natural Science Foundation of China (11271047,
11371204) and the Fundamental Research Funds for the Central University of China.

\end{document}